\documentclass[12pt,vatola]{article}

\usepackage{graphicx}

\def\C{\centerline}

\def\re#1{\par\hangindent\parindent\indent\llap{#1\enspace}\ignorespaces}

\def\no{\noindent}

\begin{document}

\C{\bf Automorphisms and Enumeration of Maps of Cayley Graphs }
\vskip 5mm \C{\bf of a Finite Group} \vskip 10mm \C{Linfan Mao}
\vskip 3mm \C{\scriptsize (Chinese Academy of Mathematics and
System Science, Beijing 100080, P.R.China)} \vskip 5mm \C{Yanpei
Liu} \C{\scriptsize (Institute of Applied Mathematics, Beijing
Jiaotong University, Beijing,100044, P.R.China) }

\vskip 10mm
\begin{minipage}{130mm}
\no{{\bf Abstract}: {\small A map is a connected topological graph
$\Gamma$ cellularly embedded in a surface. In this paper, applying
Tutte's algebraic representation of map, new ideas for enumerating
non-equivalent orientable or non-orientable maps of graph are
presented. By determining automorphisms of maps of Cayley graph
$\Gamma={\rm Cay}(G:S)$ with ${\rm Aut} \Gamma\cong G\times H$ on
locally, orientable and non-orientable surfaces, formulae for the
number of non-equivalent maps of $\Gamma$ on surfaces (orientable,
non-orientable or locally orientable) are obtained . Meanwhile,
using reseults on GRR graph for finite groups, we enumerate the
non-equivalent maps of GRR graph of symmetric groups, groups
generated by $3$ involutions and abelian groups on orientable or
non-orientable surfaces. } }

\no{{\bf Key words:} {\small  embedding, map, finite group, Cayley
graph, graphical regular representation, automorphism group,
Burnside Lemma.}}

\no{{\bf Classification:}  AMS(1991) 05C10,05C25,05C30 }
\end{minipage}

\vskip 8mm \no{\bf 1. Introduction} \vskip 5mm

\no{\it Maps} originate from the decomposition of surfaces. A
typical example in this field is the {\it Heawood map coloring
theorem}. Combinatorially, a {\it map} is a connected topological
graph $\Gamma$ cellularly embedded in a surface. Motivated by the
{\it four color problem}, the enumeration of maps on surfaces,
especially, the planar rooted maps, has been intensively
investigated by many researchers after the Tutte's pioneer work in
1962 (see [10]). By using the automorphisms of the sphere,
Liskovets gives an enumerative scheme for unrooted planar
 maps$^{[8]}$. Liskovets,  Walsh and Liskovets got many enumeration results
for {\it general planar maps, regular planar maps, Eulerian planar
maps, self-dual  planar maps} and {\it $2$-connected planar maps},
etc $^{[7]-[9]}$. Applying the well-known Burnside Lemma in
permutation groups and the Edmonds embedding scheme$^{[2]}$, Biggs
and White presented a formula for enumerating the non-equivalent
maps (also a kind of unrooted maps) of a graph on orientable
surfaces(see [1],[14],[19]), which has been successfully used for
the complete graphs, wheels and complete bipartite graphs by
determining the fix set $F_v(\alpha)$ for each vertex $v$ and
automorphism $\alpha$ of a graph$^{[14]-[15],[19]}$.

Notice that Biggs and White's formula can be only used for
orientable surfaces. For counting non-orientable maps of graphs,
new mechanism should be devised. In 1973,Tutte presented an
algebraic representation for maps on locally orientable
surface$^{([10],[17]-[18])}$. Applying the Tutte's map
representation, a general scheme for enumerating the
non-equivalent maps of a graph on surfaces can be established
(Lemma $3.1$ in section $3$), which can be used for orientable or
non-orientable surfaces. This enumeration scheme has been used to
enumerate complete maps on surfaces (orientable,non-orientable or
locally orientable) by determining all orientation-preserving
automorphisms of maps of a complete graph$^{[13]}$. In orientable
case, result is the same as in [14]. The approach of counting
orbits under the action of a permutation group is also used to
enumerate the rooted maps and non-congruent embeddings of a
graph$^{[6],[11],[16]}$. Notice that an algebraic approach for
construction non-hamiltonian cubic maps on every surface is
presented in $[12]$. The main purpose of this paper is to
enumerate the non-equivalent maps of Cayley graph $\Gamma$ of a
finite group $G$ satisfying ${\rm Aut} \Gamma =R(G)\times H\cong
G\times H$ on orientable, non-orientable or locally orientable
surfaces, where $H$ is a subgroup of ${\rm Aut} \Gamma$. For this
objective, we get all orientation-preserving automorphisms of maps
of $\Gamma$ in the Section $2$. The scheme for enumerating
non-equivalent maps of a graph is re-established in Section $3$.
Using this scheme, results for non-equivalent maps of Cayley
graphs are obtained. For concrete examples, in Section $4$, we
calculate the numbers of non-equivalent maps of GRR graphs for
symmetric groups, groups generated by $3$ involutions and abelian
groups. Terminologies and notations used in this paper are
standard. Some of them are mentioned in the following.

  All surfaces are 2-dimensional compact closed manifolds without
boundary, graphs are connected and groups are finite in the
context.

For a finite group $G$, choose a subset $S\subset G$ such that
$S^{-1}=S$ and $1_G\not\in S$, the Cayley graph $\Gamma= {\rm
Cay}(G:S)$ of $G$ with respect to $S$ is defined as follows:

$V(\Gamma) = G$;

$E(\Gamma) = \{(g,sg)| g\in G, s\in S \}$.

\no{It has been shown that $\Gamma$ is transitive, the right
regular representation $R(G)$ is a subgroup of ${\rm Aut}\Gamma$
and it is connected if and only if $G=\left< S \right>$. If there
exists a {\it Cayley set} $S$ such that ${\rm Aut}({\rm
Cay}(G:S))=R(G)\cong G$, then $G$ is called to have a {\it
graphical regular representation}, abbreviated to GRR and say
${\rm Cay}(G:S)$ is the GRR graph of the finite group $G$. Notice
that which groups have GRR are completely determined (see
$[4]-[5]$ and $[21]$ for details).}

A {\it{map}} $M = ({\cal X} _{\alpha,\beta},\cal{P})$ is defined ]
to be a permutation $\cal{P}$ acting on ${\cal X} _{\alpha,\beta}$
of a disjoint union of quadricells $Kx$ of $x\in{\cal X}$, where
$K=\{1,\alpha,\beta,\alpha\beta \}$ is the Klein group, satisfying
the following conditions:

$(i)$ \ for $\forall{x}\in {{\cal X} _{\alpha,\beta}}$, there does
not exist an integer $k$ such
  that ${\cal{P}}^{k}x = \alpha x$;

$(ii)$ \ $\alpha{\cal{P}}={\cal{P}}^{-1}\alpha$;

$(iii)$ \ the group $\Psi_{J}=\left<\alpha,\beta,\cal{P}\right>$
is transitive on ${\cal X}_{\alpha,\beta}$.

According to the condition $(ii)$, the vertices of a map are
defined to be the pairs of conjugate of ${\cal P}$ action on
${\cal X}_{\alpha,\beta}$ and edges the orbits of $K$ on ${\cal
X}_{\alpha,\beta}$. For example, $\{x,\alpha x,\beta x,\alpha\beta
x\}$ is an edge for $\forall x\in{\mathcal X}_{\alpha,\beta}$ of
M. Geometrically, any map $M$ is an embedding of a graph $\Gamma$
on a surface ( see also $[10],[17]-[18]$ ), denoted by
$M=M(\Gamma)$ and $\Gamma=\Gamma(M)$. The graph $\Gamma$ is called
the underlying graph.  If $r\in {\cal X}_{\alpha,\beta}$ is marked
beforehand, then $M$ is called a {\it{rooted map}}, denoted by
$M^{r}$.

For example, the graph $K_4$ on the tours with one face length $4$
and another $8$ shown in Fig. $1$,

\includegraphics[bb=-10 10 200 230]{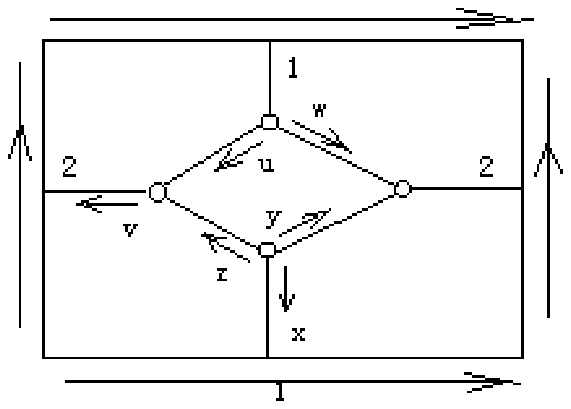}

\vskip 3mm

\C{Fig.$1$}

\vskip 3mm \no{can be algebraically represented as follows:}

{\it A map $({\cal X}_{\alpha,\beta},\cal{P})$ with ${\cal
X}_{\alpha,\beta}= \{x,y,z,u,v,w,\alpha x,\alpha y,$ $ \alpha
z,\alpha u,\alpha v,\alpha w, \beta x,\beta y,\beta z,$ $\beta u,
\beta v,\beta w,\alpha\beta x, \alpha \beta y,\alpha \beta
z,\alpha \beta u,\alpha \beta v,\alpha\beta w \}$ and}

\begin{eqnarray*}
{\cal P} &=& (x,y,z)(\alpha \beta x,u,w)(\alpha \beta z,\alpha \beta u,v)
(\alpha \beta y,\alpha \beta v,\alpha \beta w)\\
&\times& (\alpha x,\alpha z,\alpha y)(\beta x,\alpha w,\alpha u)(\beta z,\alpha v,\beta u)(\beta y,\beta w,\beta v)
\end{eqnarray*}

\no The four vertices of this map are $ \{(x,y,z), (\alpha
x,\alpha z,\alpha y)\}$, $\{(\alpha \beta x,u,w),(\beta x,\alpha
w,\alpha u)\}$, $\{(\alpha \beta z,\alpha \beta u,v),(\beta
z,\alpha v,\beta u)\}$ and $\{(\alpha \beta y,\alpha \beta
v,\alpha \beta w),(\beta y,\beta w,\beta v)\}$ and six edges are
$\{e,\alpha e,\beta e,\alpha\beta e\}$ for $\forall e\in
\{x,y,z,u,v,w\}$.

Two maps $M_{1} = ({\cal X}_{\alpha,\beta}^{1},{\cal{P}}_{1})$ and
$M_{2} = ({\cal X}_{\alpha,\beta}^{2},{\cal{P}}_{2})$ are called to be
 {\it{isomorphic}} if there exists a bijection
$\tau: {\cal X}_{\alpha,\beta}^{1} \longrightarrow {\cal
X}_{\alpha,\beta}^{2}$ such that for $\forall{x}\in{{\cal
X}_{\alpha,\beta}^{1}}$,$\tau\alpha(x)=\alpha\tau(x)$,
$\tau\beta(x)=\beta\tau(x)$ and
$\tau{\cal{P}}_{1}(x)={\cal{P}}_{2}\tau(x)$ and $\tau$ is called
an {\it{isomorphism}} between them. Similarly, two maps $M_1, M_2$
are called to be {\it equivalent} if there exists an isomorphism
$\xi$ between $M_1$ and $M_2$ such that for $\forall{x}\in{{\cal
X}_{\alpha,\beta}^{1}}$, $\tau{\cal{P}}_{1}(x)
\not={\cal{P}}_{2}^{-1}\tau(x)$. Call $\xi$ an {\it{equivalence}}
between $M_1$ and $M_2$. If $M_{1}=M_{2}=M$, then an isomorphism
or an equivalence between $M_{1}$ and $M_{2}$ is called an
automorphism or an orientation-preserving automorphism of $M$.
Certainly, an orientation-preserving automorphism of a map is an
automorphism of map preserving the orientation on this map.

All automorphisms or orientation-preserving automorphisms of a map
$M$ form groups, called {\it{automorphism group}} or {\it
orientation-preserving automorphism group} of $M$ and denoted by
${\rm Aut}M$ or ${\rm Aut_OM}$, respectively. Similarly, two
rooted maps  $M^{r}_{1}$ and $M^{r}_{2}$ are said to be
{\it{isomorphic}} if there is an isomorphism $\theta$ between them
such that $\theta(r_{1})= r_{2}$, where $r_{1}$, $r_{2}$ are the
roots of $M_{1}^{r}$ and $M_{2}^{r}$, respectively and denote the
{\it{automorphism group}} of $M^{r}$ by ${\rm AutM}^{r}$. It has
been known that ${\rm AutM}^{r}$ is the trivial group.

Now let $\Gamma$ be a connected graph. The notations ${\cal E
}^{O}(\Gamma), {\cal E}^{N}(\Gamma)$ and ${\cal E}^{L}(\Gamma)$
denote the embeddings of $\Gamma$ on the orientable surfaces,
non-orientable surfaces and locally orientable surfaces, ${\cal
M}(\Gamma)$ and $\rm Aut\Gamma$ denote the set of non-isomorphic
maps underlying a graph $\Gamma$ and its automorphism group,
respectively.

Terminologies and notations not defined here can be seen in [10]
for maps and graphs and in [1] and [20] for groups.

Notice that the {\it equivalence}  and {\it isomorphism} for maps
are two different concepts, for example, map $M = ({\cal
X}_{\alpha,\beta},\cal{P})$ is always isomorphic to its mirror map
$M^{-1} = ({\cal X} _{\alpha,\beta},{\cal P}^{-1})$, but $M_1$
must not be equivalent to its mirror $M^{-1}$. We establish an
approach for calculating non-equivalent maps underlying a graph
and concrete results in the sequel sections.

\vskip 8mm \no{\bf 2. Determining orientation-preserving
automorphisms of maps of Cayley graphs} \vskip 5mm

\no For $C=\{(x_1,x_2,\cdots,x_l),(\alpha x_l,\alpha
x_{l-1},\cdots,\alpha x_1 )\}$, the permutation $\Theta=(x_1,x_2,
\cdots,$ $x_l) (\alpha x_l,\alpha x_{l-1},\cdots,\alpha x_1 )$ is
called a pair permutation. Denote by $\{C \}$ the set
$\{x_1,x_2,\cdots,$ $x_l,\alpha x_1,\alpha x_2,\cdots,\alpha x_l
\}$ and $g\mid_{\Omega_1}$ the constraint of permutation $g$
action on $\Omega_1$ for $\Omega_1\subset\Omega$. Then we get the
following result.

\vskip 5mm
\no{{\bf Lemma 2.1} {\it Let $\Gamma$ be a connected graph. Then}}

($i$) {\it For any map $M\in {\mathcal M}(\Gamma)$, if
$\tau\in{\rm AutM}$, then $\tau\mid_{V(\Gamma)}\in{\rm
Aut\Gamma}$};

($ii$) {\it For any two maps $M_1,M_2$ underlying the graph
$\Gamma$, if $\theta$ is an isomorphism mapping $M_1$ to $M_2$,
then $\theta\mid_{V(\Gamma)}\in{\rm Aut\Gamma}$.} \vskip 3mm

{\it Proof} \ According to the Tutte's algebraic representation
for maps, we can assume that $M=({\mathcal
X}_{\alpha,\beta},{\mathcal P})$ with ${\mathcal X}=E(\Gamma)$.
For $\forall x,y\in V(M)$, we know that

$$x=\{(e_1,e_2,\cdots,e_s),(\alpha e_s,\alpha e_{s-1},\cdots,\alpha e_1)\} ;$$

$$y=\{(e^1,e^2,\cdots,e^t),(\alpha e^t,\alpha e^{t-1},\cdots,\alpha e^1)\} .$$

Now if $e=xy\in E(G)$, there must be two integers $i,j$, such that $e_i=\beta e^j=e$
or $\beta e_i=e_j=e$. Whence, we get that

($i$) if $\tau\in{\rm AutM}$, then
$V(\Gamma)=V(M)=V^{\tau}(M)=V^{\tau}(\Gamma)$ and

$$x^{\tau}=\{(\tau(e_1),\tau(e_2),\cdots,\tau(e_s)),(\alpha \tau(e_s),\alpha\tau (e_{s-1}),
\cdots,\alpha \tau(e_1))\} ;$$

$$y^{\tau}=\{(\tau(e^1),\tau(e^2),\cdots,\tau(e^t)),
(\alpha \tau(e^t),\alpha \tau(e^{t-1}),\cdots,\alpha \tau(e^1))\}.$$

\no{Therefore, }

$$e^{\tau}\in\{x^{\tau}\}\cap\beta\{y^{\tau}\} \quad {\rm or}\quad
e^{\tau}\in\beta\{x^{\tau}\}\cap\{y^{\tau}\}.$$

\no{Whence, $x^{\tau}y^{\tau}\in E(\Gamma)$ and $\tau\mid_{V(\Gamma)}\in{\rm Aut\Gamma}$.}

($ii$) \ Similarly, if $\theta : M_1 \longrightarrow M_2$ is an
isomorphism, then $\theta: V(\Gamma)=V(M_1) \longrightarrow
V(M_2)=V(\Gamma)$ and

$$e^{\theta}\in\{x^{\theta}\}\cap\beta\{y^{\theta}\} \quad {\rm or}\quad
e^{\theta}\in\beta\{x^{\theta}\}\cap\{y^{\theta}\}$$

\no{Whence we get that}

$$e^{\theta}=x^{\theta}y^{\theta}\in E(\Gamma)\quad {\rm and}\quad
\theta\mid_{V(\Gamma)}\in{\rm Aut\Gamma} .\quad\quad \natural$$

\vskip 5mm \no{{\bf Lemma 2.2} {\it For $\forall g\in{\rm AutM},
\forall x\in{{\mathcal X}_{\alpha,\beta}}$ of a map $M$,}}

($i$)\quad $|x^{\rm AutM}|=|\rm AutM|$ ;

($ii$)\quad $|x^{\prec g \succ}|= o(g)$,

\no{\it where, $o(g)$ denotes the order of $g$.} \vskip 3mm

{\it Proof} For any subgroup $H\prec {\rm AutM}$, we know that
$|H|=|x^H||H_x|$. Since $H_x \prec {\rm AutM}^x$ by definition,
where $M^x$ is the rooted map with root $x$, and ${\rm AutM}^x$ is
trivial, we know that $|H_x|=1$. Whence, $|x^H|=|H|$. Now take
$H={\rm AutM}$ or $\left< g\right>$, we get the assertions ($i$)
and ($ii$).\quad\quad $\natural$

For $\forall g\in{\rm Aut \Gamma}$, $M=({\mathcal
X}_{\alpha,\beta},{\mathcal P})\in {\mathcal M}(\Gamma)$, define
an extending action of $g$ on $M$ by
$$
g^*=g\mid^{{\mathcal X}_{\alpha,\beta}}:{\mathcal
X}_{\alpha,\beta} \longrightarrow {\mathcal X}_{\alpha,\beta},
$$
\no{such that $M^{g^*}=gMg^{-1}$ and $g\alpha =\alpha g, g\beta
=\beta g$. A permutation $p$ on set $\Omega$ is called
semi-regular if all of its orbits have the same length. Whence, an
automorphism of a map is semi-regular. The next result is followed
by Lemma $2.1$ and the definition of extending action of elements
in ${\rm Aut}\Gamma$ gives a necessary and sufficient condition
for an automorphism of a map to be an orientation-preserving
automorphism of this map.}

\vskip 5mm \no{{\bf Theorem 2.1} {\it For a connected graph
$\Gamma$, an automorphism $\xi^*$ of map $M$ is an
orientation-preserving automorphism of map underlying $\Gamma$ if
and only if there exists an element $\xi\in {\rm Aut}\Gamma$ such
that $\xi^* =\xi\mid^{{\mathcal X}_{\alpha,\beta}}$.}} \vskip 3mm

Now for a finite group $G$, let $\Gamma={\rm Cay}(G:S)$ be a connected Cayley graph respect to $S$.
Then its edge set is $\{(g,sg)|\forall g\in G, \forall s\in S \}$. For convenience, we use
$g^{sg}$ denoting an edge $(g,sg)$ in the Cayley graph ${\rm Cay}(G:S)$. Then its quadricell of
this edge can be represented by $\{g^{sg+},g^{sg-},(sg)^{g+},(sg)^{g-}\}$ and
$$
{\mathcal X}_{\alpha,\beta}(\Gamma)=\{g^{sg+}|\forall g\in
G,\forall s\in S\}\cup \{g^{sg-}|\forall g\in G,\forall s\in S\} ;
$$

$$
\alpha=\prod\limits_{g\in G,s\in S}(g^{sg+},g^{sg-});
$$

$$
\beta=\prod\limits_{g\in G,s\in S}(g^{sg+},(sg)^{g+})(g^{sg-},(sg)^{g-}).
$$

The main result of this section is the following.

\vskip 5mm \no{{\bf Theorem 2.2} {\it Let $\Gamma={\rm Cay}(G:S)$
be a connected Cayley graph with ${\rm Aut} \Gamma=R(G)\times H$.
Then for $\forall\theta\in{\rm Aut}\Gamma$, the extending action
$\theta\mid^{{\mathcal X}_{\alpha,\beta}}$ is an
orientation-preserving automorphism of a map in ${\mathcal
E}(\Gamma)$ on surfaces.}} \vskip 3mm

{\it Proof}\quad The proof is divided into two parts. First, we
prove each automorphism of the graph $\Gamma$ is semi-regular and
second, construct a stable embedding of $\Gamma$ for
$\forall\theta\in {\rm Aut \Gamma}$.

($i$) For $\forall g\in {\rm Aut} \Gamma$, since ${\rm Aut}
\Gamma=R(G)\times H$, there must exist $\gamma\in R(G), \delta\in
H$ such that $g=\gamma\delta=\delta\gamma$. Now for $\forall x\in
G$, the action of elements in $\left< g\right>$ on $x$ are as
follows.

$x^g= (x^{\delta})^{\gamma}= x^{\delta} \gamma ;$

$x^{g^2}= (x^{\delta^2})^{\gamma^2}= x^{\delta^2} \gamma^2;$

$\cdots\cdots\cdots\cdots\cdots\cdots\cdots\cdots;$

$x^{g^n}=  (x^{\delta^n})^{\gamma^n}= x^{\delta^n} \gamma^n;$

$\cdots\cdots\cdots\cdots\cdots\cdots\cdots\cdots .$

\no{Therefore, the orbit of $\left< g\right>$ acting on $x$ is}

$$
x^{\left< g \right>} = (x,x^{\delta}\gamma, x^{\delta^2}
\gamma^2,\cdots, x^{\delta^{n}} \gamma^{n},\cdots).
$$

\no{That is, for $\forall x\in G$, $|x^{\left< g
\right>}|=[o(\delta),o(\gamma)]$. Whence, $g$ is semi-regular.}

($ii$) Assume that the automorphism $\theta$ of $\Gamma$ is

$$
\theta=(a,b,\cdots,c)\cdots(g,h,\cdots,k)\cdots(x,y,\cdots,z),
$$
\no{where the length of each cycle is $\kappa=o(g)$,
$G=\{a,b,\cdots,c,\cdots,g,h,\cdots,k,\cdots,x,y,$ $\cdots,z \}$
and $S=\{s_1,s_2,\cdots,s_t\}\subset G$. Denote by
$T=\{a,\cdots,g,\cdots,x\}$ the representation set of each cycle
in $\theta$. We construct a map $M=({\mathcal
X}_{\alpha,\beta},{\mathcal P})$ underlying $\Gamma$ with}

$$
{\mathcal X}_{\alpha,\beta}(\Gamma)=\{g^{sg+}|\forall g\in
G,\forall s\in S\}\cup \{g^{sg-}|\forall g\in G,\forall s\in S\} ;
$$

$$
{\mathcal P}=\prod\limits_{g\in T}\prod\limits_{x\in C(g)}
(C_x)(\alpha C_x^{-1}\alpha^{-1}),
$$

\no{where $C(g)$ denotes the cycle containing $g$ and let $x=\theta^f (g)$, then}

$$
C_x = (\theta^f(g)^{\theta^f(s_1g+)},\theta^f(g)^{\theta^f(s_2g+)},\cdots,
\theta^f(g)^{\theta^f(s_tg+)})
$$

\no{and}

$$\alpha C_x^{-1}\alpha^{-1} = (\alpha\theta^f(g)^{\theta^f(s_tg-)},\alpha\theta^f(g)^{\theta^f(s_{t-1}g-)},\cdots,
\alpha\theta^f(g)^{\theta^f(s_1g-)}).
$$

It is clear that $M=\theta M \theta^{-1}$. According to Theorem
2.1, we know that $\theta\mid^{{\mathcal X}_{\alpha,\beta}}$ is an
orientation-preserving automorphism of map $M$.

Combining ($i$) with ($ii$), the proof is complete.\quad\quad
$\natural$

According to the {\it Rotation Embedding Scheme} for orientable
embeddings of a graph formalized by Edmonds in $[2]$, each
orientable complete map is just the case of eliminating the signs
"+, -" in our representation of maps. Whence,we get the following
result for orientable maps underlying a Cayley graph of a finite
group.

\vskip 5mm \no{{\bf Theorem 2.3}  {\it Let $\Gamma={\rm Cay}(G:S)$
be a connected Cayley graph with ${\rm Aut} \Gamma=R(G)\times H$.
Then for $\forall\theta\in {\rm Aut \Gamma}$, the extending action
$\theta\mid^{{\mathcal X}_{\alpha,\beta}}$ is an
orientation-preserving automorphism of a map in ${\mathcal
M}(\Gamma)$ on orientable surfaces.}} \vskip 3mm

Notices that a GRR graph $\Gamma$ of a finite group $G$ satisfies
${\rm Aut} \Gamma=R(G)$. Since $R(G)\cong R(G)\times\{1_{{\rm Aut}
\Gamma}\}$, by Theorems $2.2$ and $2.3$, we get all
orientation-preserving automorphisms of maps of GRR graphs of a
finite group as follows.

\vskip 5mm
\no{{\bf Corollary 2.1} {\it Let $\Gamma={\rm Cay}(G:S)$ be a connected GRR graph
 of a finite group $G$. Then
for $\forall\theta\in {\rm Aut \Gamma}$, the extending action
$\theta\mid^{{\mathcal X}_{\alpha,\beta}}$ is an
orientation-preserving automorphism of a map in ${\mathcal
M}(\Gamma)$ on locally orientable surfaces.}}

\vskip 5mm \no{{\bf Corollary 2.2}  {\it Let $\Gamma={\rm
Cay}(G:S)$ be a connected GRR of a finite group $G$. Then for
$\forall\theta\in {\rm Aut \Gamma}$, the extending action
$\theta\mid^{{\mathcal X}_{\alpha,\beta}}$ is an
orientation-preserving automorphism of a map in ${\mathcal
M}(\Gamma)$ on orientable surfaces.}}

\vskip 8mm \no{\bf 3. The enumeration of non-equivalent maps of
Cayley graphs} \vskip 5mm

\no According to Theorem 2.1, we can get a general scheme for
enumerating the non-equivalent maps of a graph $\Gamma$ on
surfaces.

\vskip 5mm \no{{\bf Lemma 3.1} {\it For any connected graph
$\Gamma$, let ${\mathcal E} \subset {\mathcal E}^{L}(\Gamma)$,
then the number $n({\mathcal E}, {\mathcal M})$ of non-equivalent
maps in ${\mathcal E}$ is}}
$$
n({\mathcal E},{\mathcal M})=\frac{1}{|\rm Aut
\Gamma|}\sum\limits_{g \in {\rm Aut \Gamma}} |\Phi (g)|,
$$
\noindent{\it where, $\Phi(g)=\{ {\mathcal{P}}| {\mathcal{P}}\in
{\mathcal E}$ and ${\mathcal{P}}^g= {\mathcal{P}} \}$.} \vskip 3mm

 {\it Proof}\quad  According to Theorem 2.1, two maps $M_1,M_2\in {\mathcal E}$
are equivalent if and only if there exists an automorphism
$g\in{\rm Aut \Gamma}$ such that $M_1^{g^*}=M_2$,
where,$g^*=g\mid^{{\mathcal X}_{\alpha,\beta}}$. Whence, all
non-equivalent maps in ${\mathcal E}$ are just the representations
of the orbits in ${\mathcal E}$ under the action of $\rm Aut
\Gamma$. By the Burnside Lemma, the number of non-equivalent maps
in ${\mathcal E}$ is
$$
n({\mathcal E},{\mathcal M})=\frac{1}{|\rm Aut
\Gamma|}\sum\limits_{g \in {\rm Aut \Gamma}} |\Phi
(g)|.\quad\quad\natural
$$

\vskip 5mm \no{{\bf Corollary 3.1} {\it The numbers of
non-equivalent maps in ${\mathcal E }^{O}(\Gamma),{\mathcal
E}^{N}(\Gamma)$ and ${\mathcal E}^{L}(\Gamma)$ are}}

$$
n({\mathcal E}^{O}(\Gamma),{\mathcal M})=\frac{1}{|\rm Aut
\Gamma|}\sum\limits_{g \in {\rm Aut \Gamma}} |\Phi^{O} (g)|;
\quad\quad (3.1)
$$

$$
n({\mathcal E}^{N}(\Gamma),{\mathcal M})=\frac{1}{|\rm Aut
\Gamma|}\sum\limits_{g \in {\rm Aut \Gamma}} |\Phi^{N} (g)|;
\quad\quad (3.2)
$$

$$
n({\mathcal E}^{L}(\Gamma),{\mathcal M})=\frac{1}{|\rm Aut
\Gamma|}\sum\limits_{g \in {\rm Aut \Gamma}} |\Phi^{L} (g)|,
\quad\quad (3.3)
$$

\noindent{\it where, $\Phi^{O}(g)=\{ {\mathcal{P}}|
{\mathcal{P}}\in {\mathcal E}^{O}(\Gamma)$ and ${\mathcal{P}}^g=
{\mathcal{P}} \}$, $\Phi^{N}(g)=\{ {\mathcal{P}}| {\mathcal{P}}\in
{\mathcal E}^{N}(\Gamma)$ and ${\mathcal{P}}^g= {\mathcal{P}} \}$,
$\Phi^{L}(g)=\{ {\mathcal{P}}| {\mathcal{P}}\in {\mathcal
E}^{L}(\Gamma)$ and ${\mathcal{P}}^g= {\mathcal{P}} \}$.}

\vskip 5mm \no{{\bf Corollary 3.2} {\it In formula (3.1)-(3.3),
$|\Phi (g)|\not=0$ if, and only if $g$ is an
orientation-preserving automorphism of map of graph $\Gamma$ on an
orientable, non-orientable or locally orientable surface.}}

\vskip 3mm The formula (3.1) is obtained by Biggs and White in
[1]. Applying Theorems $2.2-2.3$ and the formulae $(3.1)-(3.3)$,
we can enumerate the non-equivalent maps underlying a Cayley graph
$\Gamma$ of a finite group $G$ satisfying ${\rm Aut}
\Gamma=R(G)\times H$ on orientable surfaces, non-orientable
surfaces and locally orientable surfaces.

\vskip 5mm \no{{\bf Theorem 3.1} {\it Let $\Gamma={\rm Cay}(G:S)$
be a connected Cayley graph with ${\rm Aut} \Gamma=R(G)\times H$.
Then the number $n^T_{\mathcal M}(G:S)$ of non-equivalent maps
underlying $\Gamma$ on locally orientable surfaces is}}

$$
n^L_{\mathcal M} (G:S) =\frac{1}{|G||H|}\sum\limits_{\xi\in O_G}
|{\mathcal E}_{\xi}| 2^{\alpha
(S,\xi)}(|S|-1)!^{\frac{|G|}{o(\xi)}},
$$

\no{\it where $O_G$ denotes the representation set of conjugate
class of ${\rm Aut} \Gamma$ , ${\mathcal E}_{\xi}$ the conjugate
class in ${\rm Aut} \Gamma$ containing $\xi$ and}

\[
\alpha (S,\xi)=\left\{\begin{array}{ll}
\frac{|G||S|-2|G|}{2 o(\xi)},& \quad {\rm if}\quad\xi\in\Theta\\
\frac{|G||S|+2l-2|G|}{2o(\xi)},& \quad{\rm if} \quad\xi\in\Delta.
\end{array}
\right.
\]
\no{\it where, $\Theta =\{\xi| o(\xi)\equiv 1(mod 2) \vee
o(\xi)\equiv 0(mod 2), \not\exists s\in S, t\in G  \  such\ that \
s=t^{\xi^{\frac{o(\xi)}{2}}}\}$, $\Delta = \{\xi| o(\xi)\equiv
0(mod 2) \wedge \exists s_i\in S,t_i\in G, 1\leq i \leq l(\xi),
l(\xi)\equiv 0(mod \frac{o(\xi)}{2})
 \ such$ $that \ s_i=t_i^{\xi^{\frac{o(\xi)}{2}}}\}$.}

\vskip 3mm

{\it Proof}\quad Notice that $\Phi^{L}(\xi)$ is a class function
on ${\rm Aut} \Gamma$. According to Theorem 2.2 and Corollary 3.1,
we know that

\begin{eqnarray*}
n^L_{\mathcal M} (G:S) &=& \frac{1}{|\rm Aut \Gamma|}\times
\sum\limits_{\xi \in {\rm Aut} \Gamma}|\Phi^{L}(\xi)|\\
&=& \frac{1}{|G||H|}\sum\limits_{\xi\in R(G)\times H}
|\Phi^{L}(g)|.\quad\quad (3.4)
\end{eqnarray*}

Since for $\forall\xi=(\mu,\nu)\in {\rm Aut \Gamma}$, $\xi$ is
semi-regular, without loss of generality, we can assume that

$$
\xi=(a,b,\cdots,c)\cdots(g,h,\cdots,k)\cdots(x,y,\cdots,z),
$$
\no{where the length of each cycle is $o(\xi)=[o(\mu),o(\nu)]$, }

$$
{\mathcal P}=\prod\limits_{g\in T}\prod\limits_{x\in C(g)}
(C_x)(\alpha C_x^{-1}),
$$

\no{being a map underlying the graph $\Gamma$ and stable under the
action of $\xi$, $C(g)$ denotes the cycle containing $g$ and $T$
is the representation set of cycles in $\xi$. Let
$S=\{s_1,s_2,\cdots,s_k \}$ and $x=\xi^f (g)$, then}

$$
C_x = (\xi^f(g)^{\xi^f(s_1g\nu_1)},\xi^f(g)^{\xi^f(s_2g\nu_2)},\cdots,
\xi^f(g)^{\xi^f(s_tg\nu_k)}),\quad\quad (3.5)
$$

\no{with $\nu_i\in\{+,- \}, 1\leq i\leq k$.}

Notice that the quadricell adjacent to the vertex $a$ can make
$2^{|S|-1}(|S|-1)!$ pair permutations, and for each chosen pair
permutation, the pair permutations adjacent to the vertex $x,x\in
C(a)$ are uniquely determined by ($3.5$) since ${\mathcal P}$ is
stable under the action of $\xi$.

Similarly, for each given pair permutation adjacent to a vertex
$u\in T$, the pair permutations adjacent to the vertices $v,v\in
C(u)$ are also uniquely determined by ($3.5$) since ${\mathcal P}$
is stable under the action of $\xi$.

Notice that any non-orientable embedding can be obtained by
exchanging some $x$ with $\alpha x, x\in{\mathcal
X}_{\alpha,\beta} (M)$ in an orientable embedding $M$ underlying
$\Gamma$. Now for an orientable embedding $M_1$ of $\Gamma$, all
the induced embeddings by exchanging some edge's two sides and
retaining the others unchanged in $M_1$ are the same as $M_1$ by
the definition of embedding. Therefore, the number of different
stable maps under the action of $\xi$ gotten by exchanging $x$ and
$\alpha x$ in $M_1$ for $x\in U, U\subset{\mathcal X}_{\beta}$,
where ${\mathcal X}_{\beta}= \bigcup\limits_{x\in E(\Gamma)}
\{x,\beta x \}$ , is $2^{\xi (\varepsilon)-\frac{|G|}{o(\xi)}}$,
where $\xi(\varepsilon)$ is the number of orbits of $E(\Gamma)$
under the action of $\xi$, and we subtract $\frac{|G|}{o(\xi)}$
because we can choose $a^{b+},\cdots,g^{a+},\cdots, x^{a+}$ first
in our enumeration.

Since the length of each orbit under the action of $\xi$ is
$o(\xi)$ for $\forall e\in E(\Gamma)$ if $o(\xi)\equiv 1(mod 2)$
or $o(\xi)\equiv 0(mod 2)$ but there are not $s\in S, t\in G$ such
that $s=t^{\xi^{\frac{o(\xi)}{2}}}$ and is $\frac{o(\xi)}{2}$ for
each edge $t_i^{s_i t_i}$, $1\leq i\leq l(\xi)$, if $o(\xi)\equiv
0(mod 2)$ and there are $s_i \in S, t_i \in G$, $1\leq i\leq
l(\xi)$, such that $s_i= t_i^{\xi^{\frac{o(\xi)}{2}}}$ (Notice
that there must be $l\equiv 0(mod \frac{o(\xi)}{2})$ because $\xi$
is an automorphism of the graph $\Gamma$) or $o(\xi)$ for all
other edges. Whence, we get that

\[
\xi (\varepsilon)=\left\{\begin{array}{ll}
\frac{\varepsilon (\Gamma)}{o(\xi)},& \quad {\rm if}\quad\xi\in\Theta\\
\frac{\varepsilon (\Gamma) -l(\xi)}{o(\xi)}+\frac{2l(\xi)}{o(\xi)},& \quad{\rm if}\quad\xi\in\Delta .
\end{array}
\right.
\]

Now for $\forall\pi\in {\rm Aut}\Gamma$, since
$\theta=\pi\xi\pi^{-1}\in{\rm Aut} \Gamma$, we know that $\theta
(\varepsilon) = \xi (\varepsilon)$. Therefore, we get that

\[
\alpha (S,\xi)=\left\{\begin{array}{ll}
\frac{|G||S|-2|G|}{2 o(\xi)},& \quad {\rm if} \quad\xi\in\Theta\\
\frac{|G||S|+2l(\xi)-2|G|}{2o(\xi)},& \quad{\rm if} \quad\xi\in\Delta.
\end{array}
\right.
\]

\no{and}

$$|\Phi^T (\xi)|= 2^{\alpha(S,\xi)}(|S|-1)!^{\frac{|G|}{o(\xi)}}. \quad\quad (3.6)$$

Combining ($3.4$) with ($3.6$), we get that

$$
n^T_{\mathcal M} (G:S) =\frac{1}{|G||H|}\sum\limits_{\xi\in O_G}
|{\mathcal E}_{\xi}| 2^{\alpha
(S,\xi)}(|S|-1)!^{\frac{|G|}{o(\xi)}},
$$

\no{and the proof is complete.} \quad\quad $\natural$ \vskip 3mm

According to the formula ($3.1$) and Theorem $2.3$, we also get
the number $n^O_{\mathcal M} (G:S)$ of non-equivalent maps of a
Cayley graph ${\rm Cay}(G:S)$  on orientable surfaces.

\vskip 5mm \no{{\bf Theorem 3.2} {\it Let $\Gamma={\rm Cay}(G:S)$
be a Cayley graph with ${\rm Aut} \Gamma= R(G)\times H$. Then the
number $n^O_{\mathcal M}(G:S)$ of non-equivalent maps underlying
$\Gamma$ on orientable surfaces is }}

$$
n^O_{\mathcal M} (G:S)= \frac{1}{|G||H|}\sum\limits_{\xi \in O_G}
|{\mathcal E}_{\xi}|(|S|-1)!^{\frac{|G|}{o(\xi)}},
$$

\no{\it where,the means of notations ${\mathcal E}_{\xi}, O_G$ are
the same as in Theorem 3.1. } \vskip 3mm

{\it Proof} \quad By Corollary 3.1, we know that

$$
n^O_{\mathcal M} (G:S)=\frac{1}{|G||H|}\times\sum\limits_{\xi \in
{R(G)\times H}} |\Phi^{O} (\xi)|.
$$
Similar to the proof of Theorem 3.1 by applying Theorem 2.3 and
Corollary 3.1,
 we get that for $\forall \xi\in R(G)\times H$,

$$|\Phi^O (\xi)|= (|S|-1)!^{\frac{|G|}{o(\xi)}}. $$

Therefore,

$$
n^O_{\mathcal M} (G:S)= \frac{1}{|G||H|}\sum\limits_{\xi \in O_G}
|{\mathcal E}_{\xi}|(|S|-1)!^{\frac{|G|}{o(\xi)}}.\quad\quad
\natural
$$

Notice that for a given Cayley graph ${\rm Cay}(G:S)$ of a finite
group $G$, $n^O_{\mathcal M} (G:S)+$ $n^N_{\mathcal M}
(G:S)=n^L_{\mathcal M} (G:S)$. Whence, we get the number of
non-equivalent maps underlying a graph ${\rm Cay}(G:S)$ on
non-orientable surfaces.

\vskip 5mm \no{{\bf Theorem 3.3} {\it Let $\Gamma={\rm Cay}(G:S)$
be a Cayley graph with ${\rm Aut} \Gamma=R(G)\times H$. Then the
number $n^N_{\mathcal M}(G:S)$ of non-equivalent maps underlying
$\Gamma$  on non-orientable surfaces is}}
$$
n^N_{\mathcal M} (G:S) =\frac{1}{|G||H|}\sum\limits_{\xi \in O_G}
|{\mathcal E}_{\xi}|(2^{\alpha (S,\xi)}
-1)(|S|-1)!^{\frac{|G|}{o(\xi)}},
$$
\no{\it where $O_G$ denotes the representation set of conjugate
class of ${\rm Aut} \Gamma$, ${\mathcal E}_{\xi}$ the conjugate
class in ${\rm Aut} \Gamma$ containing $\xi$ and $\alpha (S,\xi)$
is the same as in Theorem 3.1.} \vskip 3mm

Since $R(G)\cong R(G)\times\{1_{{\rm Aut} \Gamma}\}$ and the
condition $s\in S, t\in G$ such that
$s=t^{\xi^{\frac{o(\xi)}{2}}}$ turns to $s=t\xi^{\frac{o(\xi)}{2}}
t^{-1}$ when ${\rm Aut}\Gamma= R(G)$, we get the number of
non-equivalent maps underlying a GRR graph of a finite group by
Theorems $3.1-3.3$ as follows.

\vskip 5mm \no{{\bf Corollary $3.3$} {\it Let $G$ be a finite
group with a GRR graph $\Gamma={\rm Cay}(G:S)$. Then the numbers
of non-equivalent maps underlying $\Gamma$ on locally orientable,
orientable and non-orientable surfaces are respective}}
$$
n^L_{\mathcal M} (G:S) =\frac{1}{|G|}\sum\limits_{g \in O_G}
|{\mathcal E}_g| 2^{\alpha_1 (S,g)}(|S|-1)!^{\frac{|G|}{o(g)}},
$$

$$
n^O_{\mathcal M} (G:S) =\frac{1}{|G|}\sum\limits_{g \in O_G}
|{\mathcal E}_g|(|S|-1)!^{\frac{|G|}{o(g)}}
$$
\no{\it and}
$$
n^N_{\mathcal M} (G:S) =\frac{1}{|G|}\sum\limits_{g \in O_G}
|{\mathcal E}_g|(2^{\alpha_1 (S,g)}
-1)(|S|-1)!^{\frac{|G|}{o(g)}},
$$
\no{\it where $O_G$ denotes the representation set of conjugate
class of $G$, ${\mathcal E}_g$ the conjugate class in $G$
containing $g$ and}

\[
\alpha_1 (S,g)=\left\{\begin{array}{ll}
\frac{|G||S|-2|G|}{2 o(g)},& \quad {\rm if} \quad g\in\Theta'\\
\frac{|G||S|+2l(g)-2|G|}{2o(g)},& \quad{\rm if} \quad g\in\Delta'.
\end{array}
\right.
\]

\no{where, $\Theta' = \{g | o(g)\equiv 1(mod 2) \vee o(g)\equiv
0(mod 2), \forall s\in S, s\not\in {\mathcal
E}_{g^{\frac{o(g)}{2}}}\}$ and $\Delta' = \{g | o(g)\equiv 0(mod
2), \exists t_i\in G, 1\leq i \leq l(g), l(g)\equiv 0(mod
\frac{o(g)}{2}) \quad {\rm such\ that}\quad t_ig^{\frac{o(g)}{2}}
t_i^{-1}\\ \in S \}.$} \vskip 3mm

Especially, if the group $G$ has odd order, then we get the
following enumeration result for maps underlying a GRR graph of
$G$.

\vskip 5mm \no{{\bf Corollary $3.4$} {\it Let $G$ be a finite
group of odd order with a GRR graph $\Gamma={\rm Cay}(G:S)$. Then
the number $n^L_{\mathcal M} (G:S)$ of non-equivalent maps of
graph $\Gamma$ on surfaces is}}
$$
n^L_{\mathcal M} (G:S) =\frac{1}{|G|}\sum\limits_{g \in O_G}
|{\mathcal E}_g|
2^{\frac{|G||S|-2|G|}{2o(g)}}(|S|-1)!^{\frac{|G|}{o(g)}}.
$$

\vskip 8mm \no{\bf 4. Examples and calculation for GRR graphs}
\vskip 5mm

Hetze and Godsil investigated GRR for solvable, non-solvable
finite groups, respectively. They proved$^{[4],[21]}$ that every
group has GRR unless it belongs to one of the following groups:

(a) abelian groups of exponent greater than 2;

(b) generalized dicyclic groups;

(c) thirteen "exceptional" groups:

(1) $Z^2_2, Z^3_2, Z^4_2$;

(2) $D_6, D_8, D_{10}$;

(3) $A_4$;

(4) $\left< a,b,c|a^2=b^2=c^2=1, abc=bca=cab \right>$;

(5) $\left< a,b|a^8=b^2=1, bab=b^5\right>$;

(6) $\left< a,b,c|a^3=c^3=b^2=1, ac=ca, (ab)^2=(cb)^2=1 \right>$;

(7) $\left< a,b,c|a^3=b^3=c^3=1, ac=ca,
bc=cb,c=a^{-1}b^{-1}ab\right>$;

(8) $Q_8\times Z_3, Q_8\times Z_4$.

Based on results in previous section, the constructions given in
$[4]-[5]$ and Corollary $3.2$, we give some calculations for the
numbers of non-equivalent maps underlying a GRR graph on surfaces
for some special groups.

\vskip 4mm
\no{{\bf Calculation 4.1}\quad\quad {\it Symmetric group $\Sigma_n$}}
\vskip 3mm

Using the notation $(\bar{k})$ denotes a partition of the integer
$n$: $(\bar{k})=k_1,k_2,\cdots,k_n)$ such that
$1k_1+2k_2+\cdots+nk_n=n$ and $lcm(\bar{k})$ the least common
multiple of the integers $1(k_1$ times), $2(k_2$ times), $\cdots$,
$n(k_n$ times), i.e, $lcm(\bar{k})=[1 (k_1 {\rm times}), 2 (k_2
{\rm times}),$ $\cdots, n (k_n {\rm times})]$ . Godsil proved
that$^{[5]}$ {\it every symmetric group $\Sigma_n$ with $n\geq 19$
has a cubic GRR with $S=\{x,y,y^{-1}\}$, where $x^2=y^3=e$}. Since
$|\Sigma_n|=n!$, we get that the numbers of non-equivalent maps
underlying a cubic GRR graph of $\Sigma_n$ are

$$
n^L_{\mathcal M} (\Sigma_n:S) =
\frac{1}{n!}\times\sum\limits_{g\in \Sigma_n} 2^{\alpha_1
(S,g)}\times 2!^{\frac{|\Sigma_n|}{o(g)}} =
\frac{1}{n!}\times\sum\limits_{g\in \Sigma_n} 2^{\alpha
(S,g)+{\frac{n!}{o(g)}}}
$$

\no{and}

\begin{eqnarray*}
n^O_{\mathcal M}
(\Sigma_n:S)&=&\frac{1}{n!}\times\sum\limits_{g\in \Sigma_n}
2!^{\frac{|\Sigma_n|}{o(g)}}\\
&=& \frac{1}{n!}\times\sum\limits_{(\bar{k})} \frac{n!}{\prod\limits_{i=1}^{n} i^{k_i}k_i!}
\times 2^{\frac{n!}{lcm(\bar{k})}}
= \sum\limits_{(\bar{k})} \frac{2^{\frac{n!}{lcm(\bar{k})}} }
{\prod\limits_{i=1}^{n} i^{k_i}k_i!},
\end{eqnarray*}

\no{and}

$$
n^N_{\mathcal M}
(\Sigma_n:S)=\frac{1}{n!}\times\sum\limits_{g\in\Sigma_n}
2^{\frac{n!}{o(g)}}(2^{\alpha_1 (S,g)}-1)
$$

For the case $n=6m+1$, we know that$^{[5]}$ $x=b_1$ if $m\equiv 1(mod 2)$
and $x=b_2$ if $m\equiv 0(mod 2)$, where

\begin{eqnarray*}
b_1 & = & (1,4)(2,n)(3,n-1)(n-6,n-3)(n-5,n-2)\\
& \times & \prod\limits_{r=1}^{m-2} (6r,6r+3)(6r+1,6r+4)(6r+2,6r+5)
\end{eqnarray*}

\no{and}

$$b_2 = b_1 (n-12,n-9).$$

Notice that $b_1\in {\mathcal E}_{[1^3 2^{3m-1}]}$ and $b_2\in
{\mathcal E}_{[1^5 2^{3m-2}]}$. We define the sets $A_1,B_1,A_2$
and $B_2$ as follows.
$$
A_1 = \{g|g\in\Sigma_n, o(g)\equiv 1(mod 2)\quad {\rm or}\quad
o(g)\equiv 0(mod 2)\quad {\rm but}\quad g^{\frac{o(g)}{2}}\not\in
{\mathcal E}_{[1^3 2^{3m-1}]} \},
$$

$$
B_1 = \{g|g\in\Sigma_n, o(g)\equiv 0(mod 2)\quad {\rm but}\quad
g^{\frac{o(g)}{2}}\not\in {\mathcal E}_{[1^3 2^{3m-1}]} \}
$$

\no{and}

$$
A_2 = \{g|g\in\Sigma_n, o(g)\equiv 1(mod 2)\quad {\rm or}\quad
o(g)\equiv 0(mod 2)\quad {\rm but}\quad g^{\frac{o(g)}{2}}\not\in
{\mathcal E}_{[1^5 2^{3m-2}]} \},
$$

$$
B_2 = \{g|g\in\Sigma_n, o(g)\equiv 0(mod 2)\quad {\rm but}\quad
g^{\frac{o(g)}{2}}\not\in {\mathcal E}_{[1^5 2^{3m-2}]} \}.
$$

For $\forall\theta\in\Sigma_n$, if $\zeta\in A_i$ or $B_i$, $i=1$
or $2$, it is clear that $\theta\zeta\theta^{-1}\in A_i$ or $B_i$.
Whence, ${\mathcal E}_{\zeta}\subset A_i$ or $B_i$. Now
calculation shows that

\[
\l(g)=\left\{\begin{array}{ll}
3!(n-3)!!,& \quad {\rm if} \quad g\in {\mathcal E}_{[1^3 2^{3m-1}]}\\
5!(n-2)!!,& \quad{\rm if} \quad g\in {\mathcal E}_{[1^5 2^{3m-2}]}  \\
0, & {\rm otherwise .}
\end{array}
\right.
\]

\no{Therefore, we have that}

\begin{eqnarray*}
n^L_{\mathcal M} (\Sigma_n:S)|_{m\equiv 1(mod 2)} &=&
\frac{\sum\limits_{g\in \Sigma_n}
2^{\alpha_1 (S,g)+{\frac{n!}{o(g)}}}}{n!}\\
&=& \frac{1}{n!}\times\sum\limits_{(\bar{k})}
\frac{n!}{\prod\limits_{i=1}^{n} i^{k_i} k_i !} \times
2^{\alpha_1 (S,(\bar{k}))+\frac{n!}{lcm (\bar{k})}}\\
&=& \sum\limits_{(\bar{k})} \frac{2^{\alpha_1
(S,(\bar{k}))+\frac{n!}{lcm (\bar{k})}}} {\prod\limits_{i=1}^{n}
i^{k_i} k_i !},
\end{eqnarray*}

\no{where,}

\[
\alpha_1 (S,(\bar{k}))=\left\{\begin{array}{ll}
\frac{n!}{2 \cdot lcm (\bar{k})},& \quad {\rm if} \quad {\mathcal E}_{(\bar{k})}\subset A_1\\
\frac{n!+12(n-3)!!}{2 \cdot lcm (\bar{k})},& \quad{\rm if} \quad
{\mathcal E}_{(\bar{k})}\subset B_1
\end{array}
\right.
\]

\no{and}

\begin{eqnarray*}
n^L_{\mathcal M} (\Sigma_n:S)|_{m\equiv 0(mod 2)} &=&
\frac{\sum\limits_{g\in \Sigma_n}
2^{\alpha'_1 (S,g)+{\frac{n!}{o(g)}}}}{n!}\\
&=& \sum\limits_{(\bar{k})} \frac{2^{\alpha'_1
(S,(\bar{k}))+\frac{n!}{lcm (\bar{k})}}} {\prod\limits_{i=1}^{n}
i^{k_i} k_i !},
\end{eqnarray*}

\no{where}

\[
\alpha'_1 (S,(\bar{k}))=\left\{\begin{array}{ll}
\frac{n!}{2\cdot lcm (\bar{k})},& \quad {\rm if} \quad {\mathcal E}_{(\bar{k})}\subset A_2\\
\frac{n!+240(n-5)!!}{2\cdot lcm (\bar{k})},& \quad{\rm if} \quad
{\mathcal E}_{(\bar{k})}\subset B_2.
\end{array}
\right.
\]

\vskip 4mm
\no{{\bf Calculation 4.2}\quad\quad {\it Group generated by 3 involutions}}
\vskip 3mm

Let $G=\left< a,b,c|a^2=b^2=c^2=e \right>$ be a finite group of
order $n$. In [5], Godsil proved that if $({\rm Aut}G)_S=e$, where
$S=\{a,b,c\}$, then $G$ has a GRR ${\rm Cay}(G:S)$. Since any
element of order $2$ must has the form $txt^{-1}, t\in G$ and
$x=a,b$ or $c$. We assume that for $\forall t\in G, tx\not= xt$,
for $x=a,b,c$. Then for $\forall g\in G$,

\[
l(g)=\left\{\begin{array}{ll}
n,& \quad {\rm if} \quad o(g)\equiv 0(mod 2)\\
0,& \quad{\rm if} \quad o(g)\equiv 1(mod 2).
\end{array}
\right.
\]

\no{Therefore, we get that}

\[
\alpha_1 (S,g)=\left\{\begin{array}{ll}
\frac{n}{2o(g)},& \quad {\rm if} \quad o(g)\equiv 1(mod 2)\\
\frac{3n}{2o(g)},& \quad{\rm if} \quad o(g)\equiv 0(mod 2),
\end{array}
\right.
\]

$$
n^L_{\mathcal M} (G:S)=\frac{\sum\limits_{o(g)\equiv 1(mod 2)}
2^{\frac{3n}{2o(g)}}+\sum\limits_{o(g)\equiv 0(mod 2)}
2^{\frac{5n}{2o(g)}}}{n},
$$

$$
n^O_{\mathcal M} (G:S)=\frac{\sum\limits_{g\in G}
2^{\frac{n}{o(g)}}}{n}
$$

\no{and}

$$
n^N_{\mathcal M} (G:S)=\frac{\sum\limits_{o(g)\equiv 1(mod 2)}
2^{\frac{n}{o(g)}}(2^{\frac{n}{2o(g)}}-1) +\sum\limits_{o(g)\equiv
0(mod 2)} 2^{\frac{n}{o(g)}}(2^{\frac{3n}{2o(g)}}-1)}{n}
$$

\vskip 4mm
\no{{\bf Calculation 4.3}\quad\quad {\it Abelian group }}
\vskip 3mm

Let $k=|S|$. It has been proved that an abelian group $G$ has GRR
if and only if $G=(Z_2)^n$ for $n=1\quad {\rm or}\quad n\geq 5$.
Now for the abelian group $G=(Z_2)^n= \left< a\right>\times\left<
b\right>\times\cdots\times\left< c\right>$, every element in $G$
has order $2$. Calculation shows that

\[
l(g)=\left\{\begin{array}{ll}
2^n,& \quad {\rm if} \quad g\in S\\
0,& \quad{\rm if} \quad g\not\in S.
\end{array}
\right.
\]

\no{Whence, we get that}

\[
\alpha_1 (S,g)=\left\{\begin{array}{ll}
(k-2)2^{n-2},& \quad {\rm if} \quad g\not\in S\\
k2^{n-2},& \quad{\rm if} \quad g\in S.
\end{array}
\right.
\]
Therefore, the numbers of non-equivalent maps underlying a GRR
graph of $(Z_2)^n$ on locally orienatble or orientable surfaces
are
\begin{eqnarray*}
n^L_{\mathcal M} ((Z_2)^n:S)&=& \frac{1}{|G|}\times
\sum\limits_{g\in (Z_2)^n} 2^{\alpha_1 (S,g)}(k-1)!^{\frac{|G|}{o(g)}}\\
&=& \frac{1}{2^n}\times (\sum\limits_{g\in S} 2^{k2^{n-2}}(k-1)!^{2^{n-1}}
+\sum\limits_{g\not\in S,g\not=e} 2^{(k-2)2^{n-2}}(k-1)!^{2^{n-1}})\\
&=& \frac{2^{k2^{n-2}}k(k-1)!^{2^{n-1}}
+(2^n-k-1)2^{(k-2)2^{n-2}}(k-1)!^{2^{n-1}}}{2^n} \\
&+&\frac{2^{(k-2)2^{n-2}}(k-1)!^{2^n}}{2^{n}},
\end{eqnarray*}
\no{and}
\begin{eqnarray*}
n^O_{\mathcal M} ((Z_2)^n:S)&=& \frac{1}{2^{n}}\times
\sum\limits_{g\in (Z_2)^n}(k-1)!^{\frac{2^n}{o(g)}}\\
&=& \frac{(k-1)!^{2^n}+(2^n-1)(k-1)!^{2^{n-1}}}{2^{n}}.
\end{eqnarray*}

\vskip 10mm
{\bf References}
\vskip 6mm
\re{[1]} N.L.Biggs and
A.T.White, {\it Permutation Groups and Combinatoric Structure},
Cambridge University Press (1979).

\re{[2]} Edmonds.J, A combinatorial representation for polyhedral surfaces,
Notices Amer. Math. Soc 7 (1960) 646.

\re{[3]} X.G.Fang, C.H. Li, J.Wang and M.Y.Xu, On cubic graphical
regular representation of finite simple groups (preprint, Peking
University,2001).

\re{[4]} C.D.Godsil, On the full automorphism group of a graph, {\it Combinatorica}, 1(1981),243-256.

\re{[5]} C.D.Godsil, The automorphism groups of some cubic Cayley graphs,{\it Europ.J.Combin}.
4(1983),25-32
\re{[6]} Jin Ho Kwak and Jaeun Lee, Enumeration of graph embeddings,
{\it Discrete Math}, 135(1994), 129-151.

\re{[7]} V.A.Liskovets, Enumeration of non-isomorphic planar maps,{\it Sel. Math. Sov}, 4;4(1985),
303-323.

\re{[8]} V.A.Liskovets, A reductive technique for Enumerating non-isomorphic planar maps,
{\it Discrete Math}, 156(1996), 197-217.

\re{[9]} V.A.Liskovets and T.R.S.Walsh, The enumeration of non-isomorphic 2-connected planar
 maps,{\it Canad.J.Math}, 35(1983), 417-435.

\re{[10]}Y.P.Liu, {\it Enumerative Theory of Maps}, Kluwer
Academic Publisher, Dordrecht / Boston / London (1999).

\re{[11]}L.F. Mao and Y.P.Liu, New automorphism groups identity of
trees, {\it Chinese advances Math}, 1(2003).

\re{[12]}L.F. Mao and Y.P.Liu, An approach for constructing
3-connected non-hamiltonian cubic map on surfaces, {\it Chinese OR
Transactions}, 4(2001),1-7.

\re{[13]}L.F. Mao and Y.P.Liu, The number of complete maps on
surfaces, {\it arXiv: math.GM/0607}.

\re{[14]} B.P.Mull,R.G.Rieper and A.T.White, Enumeration $2$-cell imbeddings of connected graphs,
{\it Proc.Amer.Math.Soc}, 103(1988), 321-330.

\re{[15]} B.P.Mull, Enumerating the orientable $2$-cell imbeddings of complete bipartite graphs,
{\it J.Graph Theory}, vol 30, 2(1999),77-90.

\re{[16]} S.Negami, Enumeration of projective -planar embeddings of graphs,
{\it Discrete Math}, 62(1986), 299-306.

\re{[17]} W.T.Tutte, {\it Graph Theory}, Combridge university Press (2001).

\re{[18]} W.T.Tutte, What is a maps? in {\it New Directions in the Theory of Graphs}
(ed.by F.Harary), Academic Press (1973), 309-325.

\re{[19]} A.T.White, {\it Graphs of Group on Surfaces- interactions and models}, Elsevier
Science B.V. (2001).

\re{[20]} Xu Mingyao et al, {\it Introduction to finite groups (II)}, Science Press,
Beijing (1999).

\re{[21]} Xu Mingyao, Automorphism groups and isomorphisms of Cayley digraphs,
{\it Discrete Math}, 182(1998), 309-319.

\end{document}